\documentclass[11pt]{amsart}
\usepackage{graphicx}
\usepackage{fullpage}
\newcommand{\C}{\mathbb C}
\newcommand{\R}{\mathbb R}
\newcommand{\N}{\mathbb N}

\newcommand{\E}{{\mathcal E}}
\newcommand{\HH}{{\mathcal H}}
\newcommand{\eps}{\varepsilon}

\newcommand{\set}[1]{\left\{#1\right\}}
\newcommand{\virg}[1]{``#1"}
\newcommand{\jo}{J_0^{<0}}
\newcommand{\joo}{J_0^{=0}}
\newtheorem{theorem}{Theorem}[section]
\newtheorem{lemma}[theorem]{Lemma}
\newtheorem{prop}[theorem]{Proposition}

\theoremstyle{definition}
\newtheorem{definition}[theorem]{Definition}

\theoremstyle{remark}

\numberwithin{equation}{section}
\begin{document}

\title[]{Existence and multiplicity of stable bound states for the nonlinear Klein-Gordon equation}%
\author{Claudio Bonanno}%
\address{Dipartimento di Matematica Applicata \virg{U. Dini}, via Buonarroti 1/c, 56127 Pisa, ITALY}%
\email{bonanno@mail.dm.unipi.it}%


\begin{abstract}
We are interested in the problem of existence of soliton-like solutions for the nonlinear Klein-Gordon equation. In particular we study some necessary and sufficient conditions on the nonlinear term to obtain solitons of a given charge. We remark that the conditions we consider can be easily verified. Moreover we show that multiplicity of solitons of the same charge is guaranteed by the ``shape'' of the nonlinear term for equations on $\R^{N}$, hence without appealing to topological or geometrical properties of the domain.
\end{abstract}
\maketitle
\section{Introduction} \label{sec:intro}

In this paper we are interested in studying the existence of \emph{solitons}, that is solutions whose energy travels as a localized
packet and exhibit orbital stability, for the nonlinear Klein-Gordon equation
\begin{equation}
\psi _{tt}-\Delta \psi + F'(|\psi|)\, \frac{\psi}{|\psi|}=0  \tag{NKG}
\label{NKG}
\end{equation}
where $\psi (t,x)\in H^{1}(\R \times \R^{N},\C)$ with $N\ge 3$,
and $F:\R^{+} \rightarrow \R$ is of class $C^{2}$.
Equation (\ref{NKG}) is the Euler-Lagrange equation of the action
functional with Lagrangian density
\begin{equation}  \label{lag}
L(\psi)=\frac{1}{2}\left( \left| \psi_t \right| ^{2}-|\nabla \psi
|^{2}\right) -F(|\psi|)
\end{equation}
which is invariant for the action of the Poincar\'{e} group and
for the gauge action of the group $S^{1}$ given by
\begin{equation} \label{gauge}
\psi \mapsto e^{i\theta }\psi \qquad \theta \in \R
\end{equation}
Noether's theorem states that any invariance for a one-parameter
group of the Lagrangian implies the existence of an integral of
motion, namely of a quantity computed on solutions which is preserved with
time. Thus equation (\ref{NKG}) has ten integrals: energy,
momentum, angular momentum and ergocenter velocity. Moreover,
another integral is given by the gauge invariance (\ref{gauge}):
charge.

A typical example of solitons for equation (\ref{NKG}) is that of orbitally stable
standing waves. A \textit{standing wave} is a finite energy
solution of the form
\begin{equation}  \label{standing-wave}
\psi(t,x)=u(x)e^{-i\omega t}, \quad u \geq 0,\ \omega \in \R^{+} .
\end{equation}
for which (\ref{NKG}) takes the form
\begin{equation} \label{static}
-\Delta u+F'(u)=\omega^{2}u
\end{equation}
Definition of orbital stability is given below in Definition \ref{orbital-stab}.

The study of solitons for equation (\ref{NKG}) has a very long
history starting with the paper of Rosen \cite{rosen68}. The first step is to study the existence of solutions of (\ref{static}). Coleman et al. \cite {coleman78} and Strauss \cite{strauss} gave the first rigorous proofs of existence of standing waves
(\ref{standing-wave}) for some particular $F's$, and later
necessary and sufficient existence conditions have been found by
Berestycki and Lions \cite{Beres-Lions}.

For what concerns the orbital stability of standing waves for (\ref{NKG}), the first results are due to Shatah \cite{shatah}. He finds a necessary and sufficient condition for orbital stability which turns out to be difficult to be verified in concrete situations. Indeed this condition is given in terms of a function $d(\omega)$ on the frequency of the standing wave (\ref{standing-wave}), and an explicit computation of this function is possible only in some particular situations, when one knows explicit properties of the standing wave. For example, this computation has been carried out in \cite{ss85} for the (\ref{NKG}) equation with the choice $F(\psi) = \frac 1 2\, |\psi|^{2} - \frac 1 p\, |\psi|^{p}$. See also
\cite{gss87} for a generalisation to general Hamiltonian fields equations of methods in \cite{shatah}. Some instability results have also been obtained in \cite{ss85} and \cite{shatah2}.

More recently in \cite{BBBM}, it has been introduced a new approach to the study of the existence of solitons for the (\ref{NKG}) equation. More attention is paid to the properties that the nonlinear term $F$ has to satisfy in order to obtain existence of solitons, independently of their frequency. What turns out to be important is instead the charge of the standing waves. In particular, solitons have been obtained as \emph{ground states} of the energy $E$ (see (\ref{y-energy})) constrained to the manifold of standing waves with fixed charge, and sufficient conditions on $F$ are given for the existence of points of minimum (see Section \ref{sec:res}).

In this paper we continue the study of \cite{BBBM}. First, fixed a charge $\sigma$, we give a necessary and sufficient condition for the existence of ground states solitons of charge $\sigma$. When this condition is not satisfied, so that the global minimum for the energy $E$ constrained to the manifold of standing waves with charge $\sigma$ is not attained, we study the existence of solitons which are \emph{bound states}, that is points of local minimum. In particular, we prove that under the same assumptions on $F$ as in \cite{BBBM}, there exists an interval of charges such that the (\ref{NKG}) equation has bound states for any charge in this interval. Moreover, we find a sufficient condition on $F$ under which bound states exist for any charge. We stress again that all the conditions we consider are very easy to be verified for a given function $F$, hence results can be easily checked for any given (\ref{NKG}) equation.

Finally, we study the existence of multiple bound states for a fixed charge. The existence of multiple positive solutions for elliptic equations as (\ref{static}) has received much attention in recent years. The multiplicity of solutions has been studied for equations on domains of $\R^{N}$ with a rich topology (\cite{benci}) or geometry (\cite{dancer}), and for equations on Riemannian manifolds (\cite{manifold}, \cite{hirano}). In this paper we obtain existence of multiple solutions $(u,\omega)$ of (\ref{static}) on $\R^{N}$ for functions $F$ with ``oscillating behaviour'' (see (\ref{comp-conn-r})), hence the multiplicity phenomenon is related only to the shape of $F$. In particular, we find multiple points of local minimum for the energy $E$ constrained to the manifold of standing waves with charge in a given interval. Notice that in our result we also vary the frequency $\omega$ of solutions of (\ref{static}). However, this approach is used in \cite{prep} to obtain multiplicity results varying only the positive function $u$.

The structure of the paper is the following. Results are described in details in Section \ref{sec:res}. All the proofs are collected in Section \ref{sec:proofs}.

\section{Statement of results} \label{sec:res}

In equation (\ref{NKG}) we consider functions $F$ of the form
\begin{equation}  \label{R-u}
F(s)= \frac 1 2 \Omega^2\, s^2 + R(s) \qquad \ s\in \R^{+}
\end{equation}
with $\Omega >0$ and $R: \R^+ \to \R$ a $C^2$ function satisfying
\begin{equation} \label{H0}
R(0)=R^{\prime}(0)=R^{\prime\prime}(0)=0 \tag{$H_{0}$}
\end{equation}
\begin{equation}  \label{H1}
R(s) \ge - \frac 1 2 \Omega^2\, s^{2} \ \ \forall \, s \in \R^{+}
\tag{$H_{1}$}
\end{equation}
\begin{equation} \label{H2}
\exists\, s_0\in \R^+ \ \ \mbox{such that} \ \ R(s_0)<0
\tag{$H_{2}$}
\end{equation}
\begin{equation} \label{H3}
|R''(s)|\leq c_{1}s^{p-2}+c_{2}s^{q-2}\ \text{ with }c_{1},c_{2}>0 \text{
for }2<p\le q<2^{\ast }  \tag{$H_{3}$}
\end{equation}
For the existence of ground states solitons, it turns out that it is important to
consider the behaviour of the nonlinear term $R$ near the origin.
In particular we say that $R$ satisfies the \emph{negativity
condition} if
\begin{equation} \label{neg-cond}
\exists\, \alpha>0 \ \mbox{and} \ 0< \eps < \frac 4N \ \ \mbox{such that} \ \ R(s) < 0 \ \ \forall\, s\in (0,\alpha) \ \ \mbox{and} \ \ \limsup\limits_{s\to 0^{+}} \frac{|R(s)|}{s^{2+\eps}} >0 \tag{NC}
\end{equation}
Instead for the existence of bound states solitons, we show that an important role is played by the \emph{zero condition}
\begin{equation} \label{zero-cond}
\exists\, s_1\in \R^+ \ \ \mbox{such that} \ \ F(s_{1})=0 \tag{ZC}
\end{equation}

We now reconsider the problem from the point of view of dynamical
systems. Let us consider the space $X= H^1(\R^N,\C)\times
L^2(\R^N,\C)$. To a function $\psi(t,x) \in H^1(\R \times \R^N,\C
)$ we can associate a vector $\mathbf{\Psi} =(\psi,\psi_t)$ which
is in $X$ for any $t\in \R$. Equation (\ref{NKG}) can then be written
in the form
\begin{equation}  \label{eq-hamil}
\frac{\partial \mathbf{\Psi}}{\partial t} = A \mathbf{\Psi} - \mathbf{
G}(\mathbf{\Psi})
\end{equation}
with
$$
A = \left(
\begin{array}{cc}
0 & 1 \\[0.2cm]
\triangle & 0
\end{array}
\right)
\qquad \qquad
\mathbf{G}(\mathbf{\Psi}) = \left(
\begin{array}{c}
0 \\[0.2cm]
F'(|\psi|)
\end{array}
\right).
$$
Together with an initial condition $\mathbf{\Psi}(0,x)$, equation
(\ref {eq-hamil}) is a Cauchy problem, and its solutions define a
flow $U$ on the phase space $X$, that is a map $U:\R \times X \to
X$ given by
$$
\mathbf{\Psi}(t,x)=U(t,\mathbf{\Psi}(0,x))
$$
Global existence for the flow $U$ is guaranteed by assumption (\ref{H1}).

As discussed in the introduction, the flow $U$ preserves energy
and charge. Energy, by definition, is the quantity which is
preserved by the time invariance of the Lagrangian (\ref{lag}). It has the
form
\begin{equation}  \label{energy}
\E(\mathbf{\Psi})=\int \left[ \frac{1}{2}\left| \psi_t \right|^{2} + \frac{1
}{2} \left| \nabla \psi \right|^{2} + F(|\psi|) \right]\, dx
\end{equation}
Charge is the quantity preserved by the gauge action
(\ref{gauge}), and is given by
\begin{equation} \label{charge}
\HH(\mathbf{\Psi}) =\mathrm{Im} \int \psi_t\, \overline{\psi }\ dx
\end{equation}

A particular class of solutions of (\ref{NKG}) (and of
(\ref{eq-hamil})) is given by standing waves, namely solutions of
the form (\ref{standing-wave}) for which the couple $(u,\omega)$
satisfies equation (\ref{static}). Moreover if $u\in
H^1(\R^N,\R^{+})$ is a solution of (\ref{static}) for any given
$\omega\in \R$, we have a standing wave solution for (\ref{NKG}).
The vector $\mathbf{\Psi}$ associated to a standing wave $\psi$
is given by $\mathbf{\Psi}=(u(x)e^{-i\omega t},-i\omega
u(x)e^{-i\omega t})$. We define the
set
$$
Y:= \left\{(u,\omega) \in H^1(\R^N,\R^+)\times \R^+\right\}
$$
which is embedded into $X$ by
$$
Y \ni (u,\omega) \mapsto (u(x)e^{-i\omega t}, -i \omega u(x)e^{-i\omega t})
\in X \quad \forall\, t \in \R
$$
and contains the standing waves.

\begin{definition} \label{orbital-stab}
A standing wave $\mathbf{\Psi}=(u(x)e^{-i\omega t},-i\omega
u(x)e^{-i\omega t})$ is \emph{orbitally stable} if the set
{
\begin{equation}  \label{gamma-orbita}
\Gamma:= \left\{ (u(x+y)\, e^{-i(\omega t -\theta)}, -i\omega
u(x+y)\, e^{-i(\omega t -\theta)}) \right\}_{y\in \R^{N},\, \theta \in \R}
\subset X
\end{equation}}
is Lyapunov stable, that is for any $\eps>0$ there exists $\delta>0$ such that if $d_{X}(\Psi,\Gamma) < \delta$ then $d_{X}( U(t,\Psi), \Gamma) < \eps$ for all $t\in \R$.
\end{definition}

Let us define on the space $H^1(\R^N,\R^+)$ the $C^{1}$
functionals
\begin{equation} \label{j0}
J_0(u) := \int_{\R^N}\ \left[ \frac 1 2\, |\nabla u(x)|^2 +
R(u(x)) \right]\, dx
\end{equation}
\begin{equation} \label{ku}
K(u) := \int_{\R^N}\ |u(x)|^2\, dx
\end{equation}
and on the space $Y$ \emph{energy} and \emph{charge} as
\begin{equation} \label{y-energy}
E(u,\omega) = J_0(u) + \frac 1 2\, (\Omega^2 + \omega^2) \, K(u)
\end{equation}
\begin{equation} \label{y-charge}
H(u,\omega) = - \omega \, K(u)
\end{equation}
Notice that $E$ and $H$ are nothing but the restriction of
(\ref{energy}) and (\ref{charge}) to the subset $Y$ of standing
waves. Without loss of generality, we are considering only the
case of negative charges.

We now give some remarks on the assumptions (\ref{H0})-(\ref{H3}). In particular we show what are the consequences of these assumptions for the properties of the functionals above. Assumptions (\ref{H0}) and (\ref{H1}) imply that the term $J_{0}(u) + \frac 1 2\, \Omega^2\, K(u)$ is non-negative and has some coercivity with respect to the $H^{1}$ norm of $u$. This is shown in Lemma \ref{p-s}. Assumption (\ref{H2}) is the crucial one for the existence of standing waves for (\ref{NKG}), see for example \cite{Beres-Lions}. In particular, it implies that there exists functions $u$ for which $J_{0}(u)<0$, which turns out to be fundamental. Finally assumption (\ref{H3}) implies that $J_{0}$ is of class $C^{1}$ on $H^{1}$ with
$$
dJ_{0}(u)\, v = \int_{\R^N}\ \left[ \nabla u(x) \cdot \nabla v(x) +
R'(u(x))\, v(x) \right]\, dx
$$
(see e.g. \cite{Beres-Lions}), and that $u \mapsto \int R'(u)$ is a compact operator on radially symmetric functions with values on $H^{-1}_{r}$ (this is shown in Lemma \ref{p-s}).

A general method for finding solutions of elliptic equations as (\ref{static}) is to look for critical points of a given functional restricted to some manifold. In this case a suitable manifold turns out to be the manifold of fixed charge, defined for $\sigma \in \R^{+}$ by
\begin{equation} \label{eq:manifold}
M_\sigma := \set{(u,\omega) \in Y \ :\ |H(u,\omega)|=\sigma}
\end{equation}
It turns out that
\begin{theorem}[\cite{BBBM}] \label{bbbm-1}
A couple $(u_{0},\omega_{0})$ is a standing wave solution of (\ref{NKG}) if and only if it is a critical point of the energy $E(u,\omega)$ constrained to the manifold $M_\sigma$
for $\sigma = |H(u_{0},\omega_{0})|$. Moreover, isolated points of local minimum of $E$ on $M_{\sigma}$ are orbitally stable standing waves.
\end{theorem}

\noindent \emph{Remark.}
The condition for the point of minimum to be isolated is essential only to have that the set of minimisers is given exactly by (\ref{gamma-orbita}). Moreover we assume that points of local minimum are always isolated, since this is certainly true for ``generic'' functionals $E$.
\vskip 0.2cm

In particular, if the energy $E$ has an isolated point of global minimum on
$M_{\sigma}$ then equation (\ref{NKG}) admits soliton solutions.
Hence in \cite{BBBM} it has been studied the problem of existence
of ground states (points of global minimum) for $E$ constrained to $M_{\sigma}$, and the main result is that
\begin{theorem}[\cite{BBBM}] \label{bbbm-2}
Assumptions (\ref{H0})-(\ref{H3}) are sufficient for the
existence of ground states for $E$ constrained to $M_{\sigma}$, hence for the existence of solitons, for charges $\sigma$ larger than a threshold $\sigma_{g}$.
\end{theorem}

Notice that when restricted to $M_{\sigma}$, energy $E$ takes the form
\begin{equation} \label{y-en-utile}
E|_{M_{\sigma}}(u,\omega) = E_{\sigma}(u) = J_{0}(u) + \frac 1 2\, \left( \Omega^{2}\, K(u) + \frac{\sigma^{2}}{K(u)} \right)
\end{equation}
Another useful functional on $Y$ is the \emph{hylomorphic ratio}
introduced in \cite{BBBM}
\begin{equation} \label{eq:hyl-ratio}
\Lambda(u,\omega) := \frac{E(u,\omega)}{|H(u,\omega)|} = \frac 1
2\, \left( \omega + \frac{\Omega^2}{\omega} \right) +
\frac{J_0(u)}{\omega\, K(u)}
\end{equation}
It is proved in \cite{BBBM} that a sufficient
condition for the existence of ground states on
$M_\sigma$, for a given $\sigma$, is that
\begin{equation} \label{eq:condition}
\inf\limits_{(u,\omega) \in M_\sigma} \ \Lambda(u,\omega) < \Omega
\end{equation}
Solitons satisfying condition (\ref{eq:condition}) have been
called \emph{hylomorphic} (see also \cite{hylo} for a discussion
of the properties of these solitons).

In this paper we first give an explicit formula for the threshold
$\sigma_{g}$ introduced in Theorem \ref{bbbm-2}.

\begin{prop} \label{prop:thresh}
Under assumptions (\ref{H0})-(\ref{H3}) on $R$, the (\ref{NKG}) equation has
hylomorphic solitons only for charges $\sigma > \sigma_{g}$, where
$$
\sigma_{g} = \inf \set{\Omega\, K(u) - \sqrt{2\, K(u)\, |J_0(u)|}
\ :\ J_0(u)<0}
$$
\end{prop}

For the existence of hylomorphic solitons for any charge, it was
proved in \cite{BBBM} that a sufficient condition is given by
(\ref{neg-cond}). We now prove that this condition is also necessary.

\begin{theorem} \label{thm:inter-segm}
Under assumptions (\ref{H0})-(\ref{H3}) on $R$, there exist hylomorphic
solitons for any given charge $\sigma$ if and only if (\ref{neg-cond}) holds.
\end{theorem}

\noindent The proof follows from Proposition \ref{prop:thresh}
showing that condition (\ref{neg-cond}) is equivalent to $\sigma_{g} = 0$.

\subsection{Existence of bound states}

When condition (\ref{neg-cond}) is not satisfied, for charges $\sigma <
\sigma_g$ there is no ground state for $E$ on the manifold
$M_\sigma$ (see Lemma \ref{no-gs} in Section \ref{sec:proofs}).
Hence for $\sigma < \sigma_g$ we look for orbitally stable bound states, that is points of local minimum for $E$ restricted to $M_{\sigma}$. It
is suggested by formula for $\sigma_g$ that an
important subset of the space $H^1(\R^N,\R^+)$ is the set
\begin{equation} \label{jo-min-0}
  \jo:=\set{u \in H^1(\R^N,\R^+)\ :\ J_0(u)<0}
\end{equation}
We first give a result on the relation between $\jo$ and the
functional $K(u)$.
\begin{lemma} \label{esiste-kbar}
There exists $\bar k >0$ such that
$$
\inf\limits_{K(u)=k}\ J_0(u) \left\{
\begin{array}{ll}
 <0 & \mbox{for}\ k> \bar k\\[0.2cm]
 =0 & \mbox{for}\ k< \bar k
 \end{array} \right.
$$
but the infimum is not attained for $k< \bar k$. Moreover condition (\ref{neg-cond}) is equivalent to $\bar k =0$.
\end{lemma}

We are now ready to state our main results. The first is that

\begin{theorem} \label{main-1}
Under assumptions (\ref{H0})-(\ref{H3}) on $R$, there exists a threshold
$\sigma_b < \sigma_g$ such that for all $\sigma \in (\sigma_b,
\infty)$ the (\ref{NKG}) equation has solitons of charge $\sigma$.
\end{theorem}

For $\sigma > \sigma_g$, the solitons are those found in
\cite{BBBM}, they are ground states, hence hylomorphic. Whereas
for $\sigma \in (\sigma_b, \sigma_g]$, the solitons are bound
states, hence not hylomorphic.

The proof of this theorem follows from the following lemmas. We
introduce the notation
\begin{equation} \label{jo-ug-0}
  \joo:=\set{u \in H^1(\R^N,\R^+)\ :\ J_0(u)=0}
\end{equation}
The first lemma shows that if we find a point which realizes the
infimum of the energy in the open set $\jo$, then we find a
critical point of the energy constrained to the manifold of fixed
charge. Moreover this critical point is a point of local minimum,
hence the associated standing wave is orbitally stable.

\begin{lemma} \label{inf-minore}
There exists a threshold $\sigma_b < \sigma_g$ such that for all
$\sigma \in (\sigma_b, \sigma_g]$ it holds
\begin{equation} \label{inf-interno}
\inf_{\jo}\, E_{\sigma}(u) < \inf_{\joo}\, E_{\sigma}(u)
\end{equation}
using (\ref{y-en-utile}).
\end{lemma}

Second, we use the classical Palais-Smale condition to prove the
existence of the point of minimum. By the classical Schwartz
symmetrization argument, for the minimization of the energy we can
restrict the problem to radial functions. Hence we only need the
Palais-Smale condition on the radial functions $H^1_r(\R^N,\R^+)$.

\begin{lemma} \label{p-s}
For all $\sigma \le \sigma_{g}$, the energy $E(u,\omega)$ restricted to the manifold of fixed charge $M_{\sigma}$ satisfies the Palais-Smale condition if $u \in \jo\cap H^1_r$.
\end{lemma}

\noindent \emph{Remark.} A key step in the proof of Lemma \ref{p-s} is to show that a Palais-Smale sequence $(u_{n}, \omega_{n})$ for $E$ (see (\ref{carica-fissa})-(\ref{quasi-crit})) admits a limit $\omega < \Omega$. In the case $\sigma > \sigma_{g}$ this can be proved for a Palais-Smale sequence which is minimising for $E$ (\cite{BBBM}). Hence using the proof of Lemma \ref{p-s}, we obtain that $E$ satisfies the local Palais-Smale condition for all $\sigma >0$. This makes the proof of Theorem \ref{main-1} contain the proof of Theorem \ref{bbbm-2}.

\vskip 0.2cm

Moreover, we study the threshold $\sigma_b$. We obtain a checkable sufficient
condition on the function $F$ for having solitons of any charge.

\begin{theorem} \label{main-2}
Under assumptions (\ref{H0})-(\ref{H3}) on $R$, if condition (\ref{zero-cond}) holds then
the (\ref{NKG}) has solitons of any charge.
\end{theorem}
Notice that if condition (\ref{neg-cond}) holds, then $\bar k =0$ and (\ref{NKG})
admits hylomorphic solitons of any charge. In case condition (\ref{neg-cond})
does not hold, we still can have solitons of any charge by
condition (\ref{zero-cond}), but they are hylomorphic only for $\sigma > \sigma_{g}$ strictly larger than zero.

\subsection{Multiplicity of bound states}

Finally we study the problem of existence of multiple orbitally
stable standing waves. This is equivalent to the existence of
multiple positive solutions for the elliptic equation
(\ref{static}). The main feature of our result is that we show
that multiplicity comes only from the shape of the nonlinear term
$F(s) = \frac 1 2\, \Omega^{2} s^{2} + R(s)$. In particular, let

\begin{equation} \label{comp-conn-r}
\set{s: R(s) <0} = C_{1} \sqcup \dots \sqcup C_{\ell} \qquad \ell \in \N \cup \set{\infty}
\end{equation}
where $C_{i}$ are disjoint open intervals. We prove that

\begin{theorem} \label{main-3}
Under assumptions (\ref{H0})-(\ref{H3}), we prove that: 

if $\ell < \infty$, there exists an interval $\Sigma \subset \R^{+}$ such that for all $\sigma \in \Sigma$, the (\ref{NKG}) equation admits at least $\ell$ orbitally stable standing waves of charge $\sigma$; 

if $\ell = \infty$, for any integer $M>0$  there exists an interval $\Sigma_{M} \subset \R^{+}$ such that for all $\sigma \in \Sigma_{M}$, the (\ref{NKG}) equation admits at least $M$ orbitally stable standing waves of charge $\sigma$.
\end{theorem}

\section{Proofs} \label{sec:proofs}

\noindent \emph{Proof of Proposition \ref{prop:thresh}.} We first
note that hylomorphic solitons $(u,\omega)$ satisfy $J_0(u)<0$.
Indeed, since
$$
\inf\limits_{\omega \in \R^+}\, \frac 1 2\, \left( \omega +
\frac{\Omega^2}{\omega} \right) = \Omega
$$
from (\ref{eq:hyl-ratio}) it follows that $\Omega + J_0(u)/(\omega
K(u)) \le \Lambda(u,\omega) < \Omega$ implies $J_0(u)<0$. Hence we
restrict to functions $u$ in $\jo$ (see (\ref{jo-min-0})).

Using (\ref{y-en-utile}), we write $\Lambda$ as a function of $u$ and $\sigma$, and from basic algebra we have that
$$
\Lambda \left( u,\sigma \right) = \frac 1 2 \,
\frac{\sigma}{K(u)} + \frac 1 2 \, \frac{\Omega^2 K(u)}{\sigma} +
\frac{J_0(u)}{\sigma} \ge \Omega
$$
if and only if
$$
\sigma \in I(u):= \R^+ \setminus \left( \Omega\, K(u) - \sqrt{2\,
K(u)\, |J_0(u)|},\ \Omega\, K(u) + \sqrt{2\, K(u)\, |J_0(u)|}
\right)
$$
Hence for a fixed $\sigma$, $\Lambda(u) \ge
\Omega$ for all $u \in \jo$ if and only if $\sigma$ is in $I(u)$
for all $u \in \jo$. It remains to identify the set $\cap_{u \in
\jo}\, I(u)$. We first show that
$$
\sup_{u\in \jo}\, \left( \Omega\, K(u) + \sqrt{2\, K(u)\, |J_0(u)|} \right) = + \infty
$$
This is obtained by the sequence
\begin{equation} \label{frittelle-negative}
u_n(x):= \left\{
\begin{array}{cl}
s_0 & \mbox{if}\ |x| \le r_n \\[0.2cm]
0 & \mbox{if}\ |x| \ge r_n +1 \\[0.2cm]
s_0 (1+r_n -|x|) & \mbox{if}\ r_n \le |x| \le r_n +1
\end{array}
\right.
\end{equation}
with $r_{n} \to \infty$, and $R(s_{0}) < 0$ as in (\ref{H2}). Indeed
$$
J_{0}(u_{n}) = \frac 1 2\, \int_{r_{n}}^{r_{n}+1}\ s_{0}^{2} r^{N-1}\, dr + \int_{0}^{r_{n}}\ R(s_{0})\, r^{N-1}\, dr + \int_{r_{n}}^{r_{n}+1}\ R(s_0 (1+r_n -r))\, r^{N-1}\, dr
$$
The first and last term are $O(r_{n}^{N-1})$, whereas the second is negative and grows as $r_{n}^{N}$. Hence for $n$ big enough $u_{n}\in \jo$. Moreover $K(u_{n})$ clearly goes to infinity as $r_{n} \to \infty$.

Hence it follows, that
$$
\cap_{u \in \jo}\, I(u) = \left(0, \inf_{u\in \jo}\, \left( \Omega\, K(u) - \sqrt{2\, K(u)\, |J_0(u)|} \right) \right]
$$
and we obtain that $\Lambda(u) \ge \Omega$ for all $u \in \jo$ if and only if
$$
\sigma \le \sigma_{g} := \inf_{u\in \jo}\, \left( \Omega\, K(u) - \sqrt{2\, K(u)\, |J_0(u)|} \right)
$$
Moreover it follows from the proofs of Lemmas \ref{inf-minore} and \ref{p-s} below that $\sigma_{g}$ is attained.
 \qed

\vskip 0.5cm

\noindent \emph{Proof of Theorem \ref{thm:inter-segm}.} We need
to show that condition (\ref{neg-cond}) is equivalent to $\sigma_{g}= 0$. The first implication is basically contained in \cite{BBBM}, Corollary 3.10. Let condition
(\ref{neg-cond}) hold for given $\alpha >0$ and $0<\eps < \frac 4 N$, and for any $\delta>0$ consider the sequence $\set{u_n}\in
H^1$ given by
\begin{equation} \label{frittelle-a-zero}
u_n(x):= \left\{
\begin{array}{cl}
s_n & \mbox{if}\ |x| \le r_n \\[0.2cm]
0 & \mbox{if}\ |x| \ge 2\, r_n\\[0.2cm]
s_n \left( 2 - \frac{|x|}{r_{n}} \right) & \mbox{if}\ r_n \le |x| \le 2\, r_n
\end{array}
\right.
\end{equation}
with $s_n \to 0^{+}$ and satisfying $R(s_{n}) \sim -c \, s_{n}^{2+\eps}$ for some $c>0$, and $r_n \to \infty$ such that $s_n^2\, \mu(B(0,2\, r_n)) \to
\delta$ as $n\to \infty$, where $\mu(B(0,r))$ denotes the Lebesgue
measure in $\R^N$ of the ball of center 0 and radius $r$. For
this sequence it holds
$$
\begin{array}{c}
0< \mu(B(0,r_n))\, s_n^2 \le K(u_n) \le \mu(B(0,2\, r_n))\, s_n^2 \\[0.3cm]
J_0(u_n) \le const\, \frac{s_n^2}{r_{n}^{2}}\, \mu(B(0,r_n)) \,  +  R(s_n)\, \mu(B(0,r_n))
\end{array}
$$
Writing
$$
J_{0}(u_{n}) \le const\, \mu(B(0,r_n))\, s_{n}^{2+\frac 4N} \ \left( 1-c\, s_{n}^{\eps-\frac 4 N}\right)
$$
for $n$ big enough it follows $u_n \in \jo$ since $\eps < \frac 4N$. Moreover
$$
\Omega\, K(u_n) - \sqrt{2\, K(u_n)\, |J_0(u_n)|} = K(u_n)\, \left(
\Omega - \sqrt{2\, \frac{|J_0(u_n)|}{K(u_n)}} \right) \
\le \delta\, \Omega
$$
for $n$ big enough, since
\begin{equation} \label{vanishing}
0 \le \frac{|J_0(u_n)|}{K(u_{n})} \le const \, \frac{1}{r_{n}^{2}} + \frac{|R(s_n)|}{s_n^2} + \frac{O(|R(s_n)|)\, \mu(B(0,2\, r_n)) }{s_n^2\, \mu(B(0,r_n))} \ \longrightarrow\ 0
\end{equation}
by (\ref{H0}). Hence $\sigma_g = 0$.

We now show that $\sigma_g =0$ implies (\ref{neg-cond}). First of all, using
(\ref{H1}), we have for $u\in \jo$
$$
0\le |J_0(u)| = \int_{\R^N}\ \left[ - \frac 1 2\, |\nabla u(x)|^2
- R(u(x)) \right]\, dx \le \int_{\R^N}\ \left[ \frac 1 2\,
\Omega^2\, u^2 - \frac 1 2\, |\nabla u(x)|^2 \right]\, dx
$$
from which
$$
\Omega\, K(u) - \sqrt{2\, K(u)\, |J_0(u)|}\ \ge\ \Omega\, K(u) -
\sqrt{ \Omega^2\, K(u)^2 - K(u)\, \int\, |\nabla u(x)|^2} \ =
$$
$$
=\ \Omega\, K(u) \, \left( 1- \sqrt{1 - \frac{\int\, |\nabla
u(x)|^2}{\Omega^2\, K(u)}} \right) \ \ge \frac 1 2\, \Omega\,
K(u)\ \frac{\int\, |\nabla u(x)|^2}{\Omega^2\, K(u)} =
\frac{\int\, |\nabla u(x)|^2}{2\, \Omega}
$$
where we have used $1-\sqrt{1-x} \ge \frac 1 2\, x$ for all $x\in
[0,1]$. It follows that if $\sigma_g =0$ then
\begin{equation} \label{grad-va-zero}
\exists\ \set{u_n}\in \jo \qquad \mbox{such that} \qquad \int\,
|\nabla u_n(x)|^2 \to 0
\end{equation}
Let us argue by contradiction that (NC) does not hold, in particular we assume that $R(s)<0$ in a given interval $(0,\alpha)$ but
$$
\limsup\limits_{s\to 0^{+}}\, \frac{|R(s)|}{s^{2+\frac 4N}} =0
$$
that is for any $\delta>0$ it holds 
\begin{equation}\label{contradd}
R(s) \ge - \delta\, s^{2+\frac 4N} \ \ \mbox{as}\ \ s\to 0^{+}
\end{equation}
First of all notice that (\ref{contradd}) together with (\ref{H3}) implies that there
exists a constant $c>0$ such that for all $u\in H^1$
\begin{equation} \label{erre-non-pos}
\int_{\R^N}\, R(u)\ \ge\ -\delta\, \int_{|u|<\alpha} |u|^{2+\frac 4 N} -c \, \int_{\R^N}\,
|u|^{2^*} \ge -\delta\, \int_{|u|<\alpha} |u|^{2+\frac 4 N} - const\, \left(\int\, |\nabla
u_n(x)|^2\right)^{\frac{2^*}{2}}
\end{equation}
where we have used the Sobolev inequality $\|u
\|_{2^*} \le const \| \nabla u \|_2$ for functions in $H^1$. We now show that under assumption (\ref{contradd}), the sequence $\set{u_{n}}$ in (\ref{grad-va-zero}) is such that $K(u_{n})$ is bounded. Indeed let $K(u_{n}) \to \infty$. Then to have $\sigma_{g}=0$ it is necessary that $|J_{0}(u_{n})| \sim \frac 12 \Omega^{2} K(u_{n})$. This follows from $|J_{0}(u_{n})| \le \frac 12 \Omega^{2} K(u_{n})$ and 
$$
\Omega\, K(u_n) - \sqrt{2\, K(u_n)\, |J_0(u_n)|} = K(u_n)\, \left(
\Omega - \sqrt{2\, \frac{|J_0(u_n)|}{K(u_n)}} \right) \to 0
$$
Moreover, $|J_{0}(u_{n})| \sim \frac 12 \Omega^{2} K(u_{n})$ implies by (\ref{grad-va-zero}) and (\ref{erre-non-pos})
$$
\frac 12 \Omega^{2} \sim \frac{-\int\, R(u_{n})}{K(u_{n)}} \le \frac{\delta\, \int_{|u|<\alpha} |u_{n}|^{2+\frac 4 N}}{K(u_{n})} + \frac{const\, \left(\int\, |\nabla
u_n(x)|^2\right)^{\frac{2^*}{2}}}{K(u_{n})} \, < \delta\, \alpha^{\frac 4 N} + o(1)
$$
which leads to a contradiction taking $\delta$ as small as necessary. 

Let $\set{u_n}$ be the sequence in (\ref{grad-va-zero}) with $K(u_{n})$ bounded, applying the inequality (see \cite{VP})
$$
\| u \|_{L^{q}} \le const\, \| u \|_{L^{2}}^{1-\frac N2 + \frac Nq}\ \| \nabla u \|_{L^{2}}^{\frac N2 - \frac Nq}
$$
which holds for $2\le q \le 2^{*}$ and $N\ge 3$, and (\ref{erre-non-pos}) it follows
$$
0 > J_0(u_n) = \frac 1 2\, \int_{\R^N}\, |\nabla u_n(x)|^2 +
\int_{\R^N}\, R(u_n)\ \ge$$
$$
\ge\ \frac 1 2\, \int_{\R^N}\, |\nabla
u_n(x)|^2 -\delta\, \int_{\R^{N}} |u|^{2+\frac 4 N} - const\, \left(\int\, |\nabla
u_n(x)|^2\right)^{\frac{2^*}{2}}
$$
$$
\ge \ \frac 1 2\, \int_{\R^N}\, |\nabla u_n(x)|^2 \left( 1 - const\, \delta K(u_{n})^{\frac 4N} -
const\, \left(\int_{\R^N}\, |\nabla
u_n(x)|^2\right)^{\frac{2^*}{2}-1} \right)
$$
The contradiction follows from the fact that the right hand side
is positive for $n$ big enough, since $K(u_{n})$ is bounded and $\delta$ is as small as necessary, and $2^* > 2$ and $\int\,
|\nabla u_n(x)|^2 \to 0$. \qed \vskip 0.5cm

\begin{lemma} \label{no-gs}
When condition (\ref{neg-cond}) is not satisfied, for charges $\sigma <
\sigma_g$ there is no ground state for $E$ on the manifold
$M_\sigma$.
\end{lemma}

\noindent \emph{Proof.} In this case $E(u,\omega) >
\sigma\,\Omega$ for all $u\in \jo$, as follows from the proof of
Proposition \ref{prop:thresh}, for all $u\in \joo$ as follows from
the proof of Lemma \ref{inf-minore}, for all $u \in H^1 \setminus
(\jo \cup \joo)$ for any charge. However it holds
$$
\inf_{M_\sigma}\ E(u,\omega) = \sigma\, \Omega
$$
which is obtained by considering the sequence
(\ref{frittelle-a-zero}) with $s_n \to 0$, $r_n \to \infty$ and such that
$K(u_{n}) \to \frac \sigma \Omega$. As in (\ref{vanishing}) it
follows
$$
\frac{J_0(u_n)}{K(u_{n})} \
\longrightarrow\ 0
$$
and
$$
\frac 1 2\, \left( \Omega^2 K(u_n) + \frac{\sigma^2}{K(u_n)}
\right) = \frac 1 2\, \left( \Omega^2 \, \left(\frac \sigma \Omega
+ o_{r_n}(1) \right) + \frac{\sigma^2}{\frac \sigma \Omega +
o_{r_n}(1)} \right) \ \longrightarrow\  \sigma\, \Omega
$$
The proof is finished by using (\ref{y-en-utile}) for
$E(u_n,\omega_n)$. \qed

\vskip 0.5cm

\noindent \emph{Proof of Lemma \ref{esiste-kbar}.} The first
statement is obtained by using the sequence $\set{u_n}$ defined in
(\ref{frittelle-negative}) to show that $\inf_{K(u)=k} J_{0}(u)$
is negative for $k$ big enough, hence there exists $\bar k$
defined as
$$
\bar k := \inf \set{k\ :\ \inf_{K(u)=k} J_{0}(u) <0}
$$
We need to show that the set where $\inf_{K(u)=k} J_{0}(u) <0$ is
an interval. This is shown by a re-scaling argument. For any $u\in
H^1$ let us define $u_\lambda(x) := u\left( \frac x \lambda
\right)$ for $\lambda >0$. Then
\begin{equation} \label{stime-riscal}
\begin{array}{c}
K(u_\lambda) = \lambda^N \, K(u) \\[0.2cm]
J_0(u_\lambda) =
\lambda^{N-2}\, \int\, \frac 1 2\, |\nabla u|^2 + \lambda^N\,
\int\, R(u)
\end{array}
\end{equation}
Hence if $J_0(u)<0$ then $J_0(u_\lambda) < (\lambda^{N-2} -
\lambda^N)\, \int\, \frac 1 2\, |\nabla u|^2$. This shows that for
all $k\ge K(u)$ there exists $\lambda \ge1$ such that
$K(u_\lambda) =k$ and $J_0(u_\lambda) < 0$.

If $k < \bar k$ then $\inf_{K(u)=k} J_{0}(u)=0$ which is obtained
by the sequence $\set{u_n}$ defined in (\ref{frittelle-a-zero})
with $s_{n}$ and $r_{n}$ such that $K(u_{n}) \to k$. This sequence is used to prove
also that if (\ref{neg-cond}) holds then $\bar k =0$.

If there exist $u$ with $K(u)=k<\bar k$ and $J_{0}(u)=0$, then by the previous re-scaling argument, we would obtain $k=\bar k$. Hence the infimum of $J_{0}$ is not attained on $K(u)<\bar k$.

Finally, if $\bar k=0$ we show that there are hylomorphic solitons
for any charge $\sigma$. Fixed $\sigma$, let $u\in \jo$ with
$K(u)= \frac \sigma \Omega$. Then $(u,\Omega) \in M_\sigma$ (see
(\ref{eq:manifold})) and
$$
\Lambda(u,\Omega) = \frac{J_0(u) + \Omega^2\, K(u)}{\sigma} =
\frac{J_0(u)}{\sigma} + \Omega < \Omega
$$
Hence $\inf_{M_\sigma} \Lambda < \Omega$ for any charge and the
proof is finished by Theorem \ref{thm:inter-segm}. \qed

\vskip 0.5cm

\noindent \emph{Proof of Lemma \ref{inf-minore}.} We use equation (\ref{y-en-utile}) for the energy $E$ and recall notation (\ref{jo-ug-0}). Using Lemma \ref{esiste-kbar}, it follows that
\begin{equation} \label{prima-su-ug-0}
\inf\limits_{\joo}\ E_{\sigma}(u)\ =\ \inf\limits_{K(u)\ge \bar k}\ \frac 1 2\, \left( \Omega^{2}\, K(u) + \frac{\sigma^{2}}{K(u)} \right)
\end{equation}
Let now study the threshold $\sigma_{g}$ defined in Proposition \ref{prop:thresh}. Introducing the notation
$$
J_k \ :=\ \left|\, \inf\limits_{K(u)=k}\ J_{0}(u) \right|
$$
we can write
\begin{equation} \label{pass-all-inf}
\sigma_{g} = \inf_{\jo}\ \left( \Omega\, K(u) - \sqrt{2\, K(u)\, |J_0(u)|} \right) = \inf_{k\ge \bar k}\ \left( \Omega\, k - \sqrt{2\, k\, J_{k}} \right)
\end{equation}
and we claim that
\begin{equation} \label{claim-thresh}
\sigma_{g}\ <\ \Omega\, \bar k
\end{equation}
The proof is finished using (\ref{claim-thresh}). Indeed, by (\ref{prima-su-ug-0}), studying the function
$$
f(k) := \Omega^{2}\, k + \frac{\sigma^{2}}{k}
$$
one obtains that
$$
f(k) \ge 2\, \sigma\, \Omega = f\left( \frac{\sigma}{\Omega} \right)
$$
and $f$ is strictly increasing for $k > \frac{\sigma}{\Omega}$. Hence, if $\sigma \le \sigma_{g}$, by (\ref{claim-thresh}) it follows that
$$
\bar k > \frac{\sigma_{g}}{\Omega} \ge \frac{\sigma}{\Omega}
$$
Using (\ref{prima-su-ug-0}) this implies that for $\sigma \le \sigma_{g}$
\begin{equation} \label{sec-su-ug-0}
\inf\limits_{\joo}\ E_{\sigma}(u)\ = \ \frac 1 2\, \left( \Omega^{2}\, \bar k + \frac{\sigma^{2}}{\bar k} \right)\ > \ \sigma\, \Omega
\end{equation}
Moreover, since $\sigma \to E_{\sigma}(u)$ is continuous for each $u$ fixed, and for $\sigma> \sigma_{g}$
$$
\inf\limits_{\jo}\ E_{\sigma}(u)\ < \sigma\, \Omega
$$
by (\ref{sec-su-ug-0}) it follows that there exists $\sigma_{b}< \sigma_{g}$ such that
$$
\inf\limits_{\jo}\ E_{\sigma}(u)\ <\ \inf\limits_{\joo}\ E_{\sigma}(u)
$$
for all $\sigma \in (\sigma_{b}, \sigma_{g}]$.

It remains to prove (\ref{claim-thresh}). This is done by using a re-scaling argument. Let us define for any $u\in \joo$ with $K(u)= \bar k$, the functions $u_{\lambda}(x):= u\left( \frac x \lambda \right)$ for $\lambda >0$. The functional $K$ and $J_{0}$ vary as in (\ref{stime-riscal}). In particular, for any $\eps >0$, we choose
$$
\lambda_{\eps} = \left( 1 + \frac{\eps}{\bar k} \right)^{\frac 1 N}
$$
from which
\begin{equation} \label{stime-riscal-eps}
\begin{array}{c}
K(u_{\lambda_{\eps}}) = \bar k + \eps \\[0.2cm]
J_0(u_{\lambda_{\eps}}) = \left( \left( 1 + \frac{\eps}{\bar k} \right) - \left( 1 + \frac{\eps}{\bar k} \right)^{1- \frac 2 N} \right) \ \int\, R(u) \ <\ 0
\end{array}
\end{equation}
Hence from (\ref{pass-all-inf}) and (\ref{stime-riscal-eps}) it follows that for any $\eps >0$
\begin{equation} \label{stima-sigma-g}
\sigma_{g}\ \le\ \Omega\, (\bar k + \eps) - \sqrt{2\, (\bar k + \eps)\, J_{\bar k + \eps}} \le \Omega\, (\bar k + \eps) - \sqrt{2\, (\bar k + \eps)\, |J_0(u_{\lambda_{\eps}})|}
\end{equation}
Moreover, as $\eps \to 0$, we have
$$
|J_0(u_{\lambda_{\eps}})| \sim \ \frac 2 N\, \frac{\eps}{\bar k} \ \left| \int\, R(u) \right|
$$
hence
$$
\Omega\, (\bar k + \eps) - \sqrt{2\, (\bar k + \eps)\, |J_0(u_{\lambda_{\eps}})|}\ \sim\ \Omega\, \bar k + \sqrt{\eps} \left( \Omega\, \sqrt{\eps} - \frac{2}{\sqrt{N}}\, \sqrt{\left| \int\, R(u) \right|} \right)
$$
and the right hand side is less than $\Omega\, \bar k$ for $\eps$ small enough. Hence by (\ref{stima-sigma-g}), we obtain the claim (\ref{claim-thresh}). \qed

\vskip 0.5cm

\noindent \emph{Proof of Lemma \ref{p-s}.} The first part of the proof follows \cite{BBBM}, Lemma 2.7. We repeat it here for sake of completeness. For any given charge $\sigma \le \sigma_{g}$, we have to prove that if $\set{(u_{n},\omega_{n})}$ is a sequence such that
\begin{eqnarray}
& (u_{n}, \omega_{n}) \in M_{\sigma} = \set{|H(u, \omega)| = \sigma} \qquad \forall\, n \label{carica-fissa} \\[0.2cm]
& u_{n} \in \jo \cap H^{1}_{r} \qquad \forall\, n \label{u-n-jo-radiali} \\[0.2cm]
& E(u_{n}, \omega_{n}) \mbox{ is bounded} \label{en-limitata} \\[0.2cm]
& dE\big|_{M_{\sigma}} (u_{n}, \omega_{n})\ \longrightarrow_{n\to \infty}\ 0  \label{quasi-crit}
\end{eqnarray}
then, up to a choice of a sub-sequence,
\begin{equation} \label{finale-p-s}
(u_{n},\omega_{n}) \ \stackrel{H^{1}_{r} \times \R}{\longrightarrow}_{n\to \infty}\ (u,\omega)
\end{equation}
First of all from (\ref{en-limitata}) it follows that there exists $(u,\omega)$, with $u\in H^{1}_{r}$ such that
\begin{equation} \label{limiti-deb-for}
u_{n}\  \stackrel{H^{1}_{r}}{\rightharpoonup} \ u \quad \mbox{and} \quad \omega_{n} \longrightarrow \omega < \Omega
\end{equation}
This follows from the property of coercivity of the energy. We need to show that $(u_{n},\omega_{n})$ is bounded in $H^{1} \times \R$. By (\ref{y-energy}), (\ref{j0}) and (\ref{ku}) we can write
\begin{equation} \label{en-succ}
E(u_{n}, \omega_{n}) = \int\ \left( \frac 1 2\, |\nabla u_{n}|^{2} + F(u_{n}) \right)\, dx\ +\ \frac{\sigma\, \omega_{n}}{2}
\end{equation}
where $F$ is defined in (\ref{R-u}) and by (\ref{H1}) is non-negative. Hence $\set{\omega_{n}}$ and $\set{\int\, |\nabla u_{n}|^{2}}$ are bounded. It remains to show that $\set{K(u_{n})}$ is bounded too. Let us assume on the contrary that $(u_n,\omega_n)$
satisfies (\ref{carica-fissa})-(\ref{en-limitata}) and
\begin{equation*}
\int |u_n|^2\, dx \rightarrow \infty.
\end{equation*}
From (\ref{R-u}) and (\ref{H0}) it follows that
\begin{equation}  \label{vicino-zero}
\exists\, \delta>0 \ \exists\, c_1>0, \text{ such that } F(s) \geq c_1
s^2 \text { for } 0\leq s \leq \delta.
\end{equation}
and one can choose $c_1=\min\limits_{0\leq s \leq \delta} \frac{F(s)}{\frac 1 2\, \Omega^{2}\, s^2}$.
By (\ref{en-succ}) and (\ref{en-limitata}), the term $\int F(u_n)\,dx$ is bounded and by (\ref{H0}) and (\ref{vicino-zero})
\begin{equation}  \label{for}
\int F(u_n)\, dx \geq \int_{0 \leq u_n\leq \delta} F(u_n)\, dx \geq c_1\,
\int_{0 \leq u_n \leq \delta} u_n^2\, dx .
\end{equation}
On the other hand
\begin{equation*}
\int_{0\leq u_n\leq \delta} u_n^2\, dx +\int_{u_n\geq \delta} u_n^2\, dx
\rightarrow \infty,
\end{equation*}
thus we have, by equation (\ref{for})
\begin{equation*}
\int_{u_n\geq \delta} u_n^2\, dx \rightarrow \infty.
\end{equation*}
This drives to a contradiction because
\begin{equation*}
\frac{1}{\delta^{2^*-2}} \ \int_{u_n\geq \delta} u_n^{2^*}dx \geq
\int_{u_n\geq \delta} u_n^2\, dx
\end{equation*}
and by the Sobolev embedding theorem
\begin{equation*}
\int_{u_n\geq \delta} u_n^{2^*}dx \leq \int u_n^{2^*}dx \leq const \int |\nabla
u_n|^2dx< const.
\end{equation*}
Now (\ref{limiti-deb-for}) for $\set{u_{n}}$ follows from classical arguments. We recall that the spaces $L^{p}_{r}$ for $2< p<2^{*}$ are compactly embedded in $H^{1}_{r}$, hence up to a sub-sequence,
\begin{equation} \label{conv-l-p}
u_{n}\ \stackrel{L^{p}_{r}}{\longrightarrow}\ u \ \mbox{ for }\ 2< p<2^{*}
\end{equation}
It remains to show that $\omega < \Omega$. Since $\sigma \le \sigma_{g}$, it follows from (\ref{carica-fissa}) and (\ref{u-n-jo-radiali}) that for all $n$
$$
\omega_{n} = \frac{\sigma}{K(u_{n})} \le \frac{\inf_{\jo} \left( \Omega\, K(u) - \sqrt{2\, K(u)\, |J_0(u)|}\right)}{K(u_{n})} < \Omega
$$
Hence $\omega \le \Omega$. Moreover $\omega = \Omega$ and $K(u_{n}) \ge \bar k$, as deduced from Lemma \ref{esiste-kbar}, imply that
$$
\sigma_{g} \ge \sigma = \omega_{n}\, K(u_{n}) = \liminf_{n}\, (\omega_{n}\, K(u_{n})) \ge \Omega\, \bar k
$$
which is in contradiction with (\ref{claim-thresh}). Hence (\ref{limiti-deb-for}) is proved. In particular, it follows that $\set{\omega_{n}}$ and $\set{K(u_{n})}$ are bounded away from zero.

To finish the proof we write assumption (\ref{quasi-crit}) as
\begin{eqnarray}
& < dJ_{0}(u_{n}),v > + \left( \Omega^{2} + \omega_{n}^{2} - 2 \lambda_{n}\, \omega_{n}\right)\, \int\, u_{n}\, v\ =\ <\eps_{n}, v>\ \longrightarrow\ 0 \label{der-en} \\[0.3cm]
& \left( \omega_{n} - \lambda_{n} \right)\, \alpha \, K(u_{n})\ =\ <\eta_{n}, \alpha> \ \longrightarrow\ 0 \label{der-car}
\end{eqnarray}
for all $v \in H^{1}$ and $\alpha \in \R$, where $\lambda_{n}$ is a sequence of Lagrange multipliers, and $\eps_{n} \in H^{-1}$, $\eta_{n} \in \R$. Since $\set{K(u_{n})}$ is bounded away from zero, from (\ref{der-car}) we obtain $\lambda_{n} = \omega_{n} + \delta_{n}$, with $\delta_{n} \to 0$. Hence substituting in (\ref{der-en}), we have
\begin{equation} \label{der-en-utile}
< dJ_{0}(u_{n}),v > + \left( \Omega^{2} - \omega^{2} \right) \, \int\, u_{n}\, v\ = \ <\eps_{n}, v> - 2 \delta_{n}\, \omega_{n} \, \int\, u_{n}\, v + \left( \omega_{n}^{2} - \omega^{2} \right) \, \int\, u_{n}\, v \ \longrightarrow\ 0
\end{equation}
for all $v \in H^{1}$. Writing (\ref{der-en-utile}) for two functions $u_{n}$ and $u_{m}$ and subtracting yields
\begin{equation} \label{der-duale-zero}
dJ_{0}(u_{n}) - dJ_{0}(u_{m}) + \left( \Omega^{2} - \omega^{2} \right) \, (u_{n}- u_{m}) \ \stackrel{H^{-1}}{\longrightarrow}_{n,m \to \infty}\ 0
\end{equation}
Hence, since $(u_{n} - u_{m})$ is bounded in $H^{1}$, applying (\ref{der-duale-zero}) to $(u_{n} - u_{m})$, we have
\begin{equation} \label{serve-p-s}
\int_{\R^{N}}\ \left( | \nabla u_{n} - \nabla u_{m} |^{2} + (R'(u_{n}) - R'(u_{m}))\, (u_{n}-u_{m}) + \left( \Omega^{2} - \omega^{2} \right) \, |u_{n}- u_{m}|^{2} \right)\ dx \longrightarrow_{n,m \to \infty}\ 0
\end{equation}
Writing
$$
\left| \int\ (R'(u_{n}) - R'(u_{m}))\, (u_{n}-u_{m}) \right| \ \le\ \int\ |R''(u_{m} + \theta\, (u_{n}- u_{m}))|\, |u_{n}- u_{m}|^{2}
$$
for some $\theta \in (0,1)$, and using (\ref{H3}) and the inequality
$$
\int\ \left( |u_{m}| + \theta\, |u_{n}- u_{m}| \right)^{p-2}\, |u_{n}- u_{m}|^{2} \ \le \ 2^{p-2}\, \int\  \left( |u_{m}|^{p-2} + \theta\, |u_{n}- u_{m}|^{p-2} \right)\, |u_{n}- u_{m}|^{2} \ \le
$$
$$
\le \ 2^{p-2}\, \left( \int\  |u_{m}|^{p} \right)^{1-\frac 2 p}\ \left( \int\ |u_{n}- u_{m}|^{p} \right)^{\frac 2 p} + 2^{p-2}\, \int\ |u_{n}- u_{m}|^{p}
$$
from (\ref{conv-l-p}) it follows that
$$
\int_{\R^{N}}\ (R'(u_{n}) - R'(u_{m}))\, (u_{n}-u_{m}) \ dx\ \longrightarrow_{n,m \to \infty}\ 0
$$
which in particular implies compactness for $u\mapsto \int R'(u)$ on $H^{1}_{r}$. Hence from (\ref{serve-p-s}) we obtain
$$
\| u_{n} - u_{m} \|_{H^{1}} \ \longrightarrow_{n,m \to \infty}\ 0
$$
since $\left( \Omega^{2} - \omega^{2} \right) >0$. Hence $\set{u_{n}}$ is a Cauchy sequence in $H^{1}_{r}$, and it follows that it has a convergent subsequence and (\ref{finale-p-s}) is proved. \qed

\vskip 0.5cm

\noindent \emph{Proof of Theorem \ref{main-1}.}  For $\sigma > \sigma_{g}$ there exist solitons with charge $\sigma$ by Proposition \ref{prop:thresh} and the main result of \cite{BBBM}. Let $\sigma \in (\sigma_{b}, \sigma_{g}]$, where $\sigma_{b}$ is the threshold found in Lemma \ref{inf-minore}. Let
$$
m := \inf_{u\in \jo,\ (u,\omega)\in M_{\sigma}}\ E(u,\omega)
$$
By the Schwartz spherical rearrangement of functions in $H^{1}$, the value $m$ is not changed if we restrict ourselves to functions $u \in H^{1}_{r}$. Moreover, by the Palais principle of symmetric criticality \cite{palais}, since the energy $E$ is invariant under rotations of $\R^{N}$, it follows that critical points for $E$ restricted to $M_{\sigma}\cap \set{u \in H^{1}_{r}}$ are also critical points for $E$ restricted to $M_{\sigma}$.

Since $J_{0}$ is continuous, the set $\jo$ is open in $H^{1}_{r}$, and by Lemma \ref{inf-minore} it follows that the infimum value $m$ is realised by a minimisation sequence in the interior of $\jo$. Hence we can apply the classical Ekeland principle to find a sequence $(u_{n}, \omega_{n}) \in M_{\sigma}$, with $u_{n}\in \jo$, satisfying (\ref{carica-fissa})-(\ref{quasi-crit}). Since by Lemma \ref{p-s}, the energy $E$ satisfies the Palais-Smale condition, it follows that there exists $(u,\omega) \in M_{\sigma}$ such that $E(u,\omega) = m$. Hence $(u,\omega)$ is a critical point of the energy $E$ restricted to $M_{\sigma}$, in particular it is a point of local minimum. It follows that $(u,\omega)$ is an orbitally stable standing wave with charge $\sigma$. \qed

\vskip 0.5cm

\begin{lemma} \label{prop:thresh-2}
Let condition (\ref{neg-cond}) not to hold, so that $\bar k >0$. If
\begin{equation} \label{small-omega}
\Omega^{2} < \sup\limits_{\jo}\ \frac{2|J_{0}(u)|}{K(u)-\bar k}
\end{equation}
then $\sigma_b=0$.
\end{lemma}

\noindent \emph{Proof.} In the proof of Theorem \ref{main-1}, the value of $\sigma_{b}$ is determined by the values of the charge for which (\ref{inf-interno}) holds. Hence we only need to show that (\ref{small-omega}) implies (\ref{inf-interno}) for all $\sigma \in (0,\sigma_{g}]$. Using (\ref{sec-su-ug-0}) we only need to show that (\ref{small-omega}) implies
\begin{equation} \label{tutte-car}
\inf_{\jo}\ E_{\sigma}(u) \ <\ \frac 1 2\, \left( \Omega^{2}\, \bar k + \frac{\sigma^{2}}{\bar k} \right) \qquad \forall\, \sigma \in (0,\sigma_{g}]
\end{equation}
By (\ref{small-omega}) for any $\eps >0$ there exists $u_{\eps} \in \jo$ such that
$$
\Omega^{2} \ < \ \frac{2\, |J_{0}(u_{\eps})|}{K(u_{\eps}) - \bar k}
$$
hence
$$
J_{0}(u_{\eps}) + \frac 1 2\, \Omega^{2}\, K(u_{\eps})\ <\ \frac 1 2\, \Omega^{2}\, \bar k
$$
and since $K(u_{\eps}) \ge \bar k$ by Lemma \ref{esiste-kbar} we find
$$
E_{\sigma}(u_{\eps}) = J_{0}(u_{\eps}) + \frac 1 2\, \Omega^{2}\, K(u_{\eps}) + \frac 1 2\, \frac{\sigma^{2}}{K(u_{\eps})} \ <\ \frac 1 2\, \Omega^{2}\, \bar k + \frac 1 2\, \frac{\sigma^{2}}{\bar k}
$$
for all $\sigma >0$. Then (\ref{tutte-car}) follows. \qed

\vskip 0.5cm

\noindent \emph{Proof of Theorem \ref{main-2}.} We need to show that condition (\ref{zero-cond}) implies (\ref{small-omega}). Let $s_{1}$ be defined as in (\ref{zero-cond}) and consider the sequence $\set{u_{n}}$ defined as in (\ref{frittelle-negative})
$$
u_n(x):= \left\{
\begin{array}{cl}
s_1 & \mbox{if}\ |x| \le r_n \\[0.2cm]
0 & \mbox{if}\ |x| \ge r_n +1 \\[0.2cm]
s_1 (1+r_n -|x|) & \mbox{if}\ r_n \le |x| \le r_n +1
\end{array}
\right.
$$
with $r_{n} \to \infty$, and $R(s_{1}) = - \frac 1 2\, \Omega^{2}\, s_{1}^{2}$. Then
$$
J_{0}(u_{n}) = \frac 1 2\, \int_{r_{n}}^{r_{n}+1}\ s_{1}^{2} r^{N-1}\, dr + \int_{0}^{r_{n}}\ R(s_{1})\, r^{N-1}\, dr + \int_{r_{n}}^{r_{n}+1}\ R(s_1 (1+r_n -r))\, r^{N-1}\, dr
$$
is negative for $r_{n}$ big enough and
$$
K(u_{n}) = s_{1}^{2}\, \mu(B(0,r_{n})) + o(r_{n}^{N})
$$
This implies
$$
\lim_{n\to \infty}\ \frac{2\, |J_{0}(u_{n})|}{K(u_{n})} \ =\ \Omega^{2}
$$
Hence given $\bar k$, for $n$ big enough
$$
\Omega^{2}\ < \ \frac{2\, |J_{0}(u_{n})|}{K(u_{n})-\bar k}
$$
and  (\ref{small-omega}) follows. \qed

\subsection{Proof of Theorem \ref{main-3}.} By the same argument as in the proof of Theorem \ref{main-1}, using the Schwartz spherical rearrangement of functions in $H^{1}$ and the Palais principle of symmetric criticality, we can restrict ourselves to $u \in
H^{1}_{r}$. Hence we look for points of local minimum of $E(u,\omega)$ on $M_\sigma$ and with $u \in \jo \cap H^{1}_{r}$.

The restriction to radially symmetric functions is fundamental. We recall that for functions in $H^{1}_{r}(\R^{N})$ there exist positive constants $\beta$ and $ \gamma$ only depending on $N$ such that
\begin{equation} \label{strauss}
|u_{n}(x)|\ \le \ \gamma\, \frac{\| u
\|_{H^{1}}}{|x|^{\frac{N-1}{2}}} \qquad \mbox{for} \qquad |x| \ge
\beta
\end{equation}
For a proof of this property see e.g. \cite{Beres-Lions}.

Given an interval $I\subset \R^{+}$ and a function $u\in H^{1}_{r}(\R^{N},\R^{+})$, we introduce the notation
\begin{equation} \label{funz-caratt}
\chi_{u,I} (x) := \chi_{\set{u(x)\in I}} (x)
\end{equation}
for the characteristic function of the set $\set{x:u(x)\in I}$, and
\begin{equation} \label{u-ristrett}
u_{I}(x) := \left\{
\begin{array}{ll}
u(x) & u(x)\in I \\
0 & u(x) \not\in I
\end{array} \right. \qquad u_{I} \in H^{1}_{r}\left( \set{x:u(x)\in I} \right)
\end{equation}

\begin{lemma} \label{conv-fort-caratt}
For any interval $I=(a,b)$ with $a>0$, if a sequence $\set{u_n}$ converges to $u$ in the $H^1_r$ norm then $\set{\chi_{u_n,I}(x)}$ converges to $\set{\chi_{u,I}(x)}$ in the
$L^1$ norm. Hence in particular the symmetric difference $\set{u_{n}(x)\in I} \bigtriangleup \set{u(x)\in I}$ has vanishing measure.
\end{lemma}

\noindent \emph{Proof.} Up to the choice of a sub-sequence, the sequence $\set{u_n}$ converges to $u$ for almost all $x$, hence $\set{\chi_{u_n,I}(x)}$ converges to $\set{\chi_{u,I}(x)}$ for almost all $x$. The proof is finished by the Lebesgue theorem of dominated convergence, since by (\ref{strauss}) there exists $\bar \beta \ge \beta$ such that if $|x|\ge \bar \beta$ then $|u_{n}(x)| < a$ for all $n$. Hence for all $n$ it follows $\chi_{u_n,I}(x) \le \chi_{\set{|x|\le \bar \beta}} (x) \in L^{1}$. \qed

\begin{lemma} \label{cont-jo-ristr}
For any interval $I=(a,b)$ with $a>0$, the function
$u \mapsto J_0(u_{I})$ is continuous in the $H^1_r$ norm.
\end{lemma}

\noindent \emph{Proof.} We recall that for $u\in H^{1}_{r}$ using notation (\ref{u-ristrett})
$$
J_{0}(u_{I}) = \int_{\set{u(x)\in I}}\ \left( \frac 1 2\, |\nabla u(x)|^{2} + R(u(x)) \right)\, dx
$$
Let $\set{u_n}$ be a sequence converging to $u$ in the $H^1_r$ norm. Then using notation (\ref{funz-caratt})
$$
\left| \int_{\set{u_{n}(x)\in I}}\ |\nabla u_{n}|^{2} - \int_{\set{u(x)\in I}}\ |\nabla u|^{2} \right| = \left| \int_{\R^{N}}\ \left( \chi_{u_n,I} |\nabla u_{n}|^{2} - \chi_{u,I} |\nabla u|^{2} \right) \right| \le
$$
$$
\le \int_{\R^{N}}\ |\chi_{u_n,I}|\, \left| |\nabla u_{n}|^{2} -  |\nabla u|^{2} \right| + \int_{\R^{N}}\ \left| \chi_{u_n,I} - \chi_{u,I} \right|\,  |\nabla u|^{2}
$$
The first term in the right-hand side is vanishing since $\set{u_n}$ converges to $u$ in the $H^1_r$ norm and $|\chi_{u_n,I}| \le 1$. The second term is vanishing by Lemma \ref{conv-fort-caratt}, which implies that the symmetric difference $\set{u_{n}(x)\in I} \bigtriangleup \set{u(x)\in I}$ has vanishing measure as $n\to \infty$, and by the absolute continuity of the integral, being $|\nabla u|^{2}$ in $L^{1}$. The same argument applies to the part with $R(u)$, by using the continuity of $J_{0}$. \qed

\vskip 0.5cm

Let us now recall notation (\ref{comp-conn-r}) for the disjoint intervals $C_{i}$ of the set $\set{s:R(s)<0}$, and write
\begin{equation} \label{estr-int}
C_{i} = (\xi_{i}, \eta_{i}) \qquad i=1,\dots,\ell
\end{equation}
where
$$
0 \le \xi_{1} < \eta_{1} < \xi_{2} < \dots < \xi_{i} < \eta_{i} < \xi_{i+1} < \dots < \eta_{\ell} \le \infty
$$
We start with the case $\ell$ finite. We recall from \cite{BBBM} the following result
\begin{lemma}[\cite{BBBM}] \label{principio-massimo}
Let $R$ satisfy (\ref{H0})-(\ref{H3}) and assume that there exist $\bar s>0$ such that $R'(s)\ge 0$ for $s\ge \bar s$, then for any $0<\omega< \Omega$ all solutions $u$ of (\ref{static}) satisfy
$$
\| u(x) \|_{L^{\infty}(\R^{N})} \le \bar s
$$
\end{lemma}

\begin{lemma} \label{primo-minimo}
There exists a threshold $\sigma_{g,1}$ such that, for all $\sigma > \sigma_{g,1}$, there exists a point of local minimum for the energy $E_{\sigma}(u)$ with $\| u(x) \|_{L^{\infty}} < \eta_{1}$.
\end{lemma}

\noindent \emph{Proof.} Let us consider the modified nonlinear term
\begin{equation} \label{tilde-r}
\tilde R(s) = \left\{
\begin{array}{ll}
R(s) & s\le \eta_{1} \\[0.2cm]
f(s) & s\ge \eta_{1}
\end{array} \right.
\end{equation}
for any function $f(s)$ for which $f'(s) \ge 0$ and such that $\tilde R$ is of class $C^{2}$ and satisfies (\ref{H0})-(\ref{H3}). Then $\tilde R$ satisfies the assumptions of Lemma \ref{principio-massimo} with $\bar s = \eta_{1}$.

We first apply Proposition \ref{prop:thresh} (or Theorem \ref{bbbm-2}) to $\tilde R$ as defined in (\ref{tilde-r}). Hence we find a threshold $\sigma_{g,1}$ such that the (\ref{NKG}) has hylomorphic solitons $(\tilde u, \tilde \omega)$ with charge $\sigma > \sigma_{g,1}$. Moreover, recalling (\ref{limiti-deb-for}) it holds $\tilde \omega < \Omega$. We can now apply Lemma \ref{principio-massimo} and obtain $\| \tilde u(x) \|_{L^{\infty}} < \eta_{1}$. 

In particular $\tilde u$ is a point of local minimum for
$$
\tilde E_{\sigma}(u) = \int\ \left( \frac 1 2\, |\nabla u|^{2} + \tilde R(u) \right)\, dx + \frac 1 2\, \left( \Omega^{2}\, K(u) + \frac{\sigma^{2}}{K(u)} \right)
$$
We now argue that $\tilde u$ is a point of local minimum also for the energy $E_{\sigma}$ with nonlinear term $R$, hence we choose $u = \tilde u$. 

Recall from Lemma \ref{conv-fort-caratt} that functions $v$ close to $\tilde u$ in the $H^{1}$ norm satisfy
$$
m(\set{v(x) \in (\eta_{1}, \infty)}) \longrightarrow 0 \qquad \mbox{ as } \ \| v - \tilde u \|_{H^{1}} \to 0
$$
This implies in particular that for all $i= 1, \dots, \ell$
$$
\int_{\set{v(x) \in (\eta_{i}, \eta_{i+1})}}\ |v|^{2} \ \le \ 2\, \left( \| v - \tilde u \|^{2}_{H^{1}} +  \int_{\set{v(x) \in (\eta_{i}, \eta_{i+1})}}\ |u|^{2} \right) \ \ \longrightarrow 0 \qquad \mbox{ as } \ \| v - \tilde u \|_{H^{1}} \to 0
$$
whence, by notation (\ref{u-ristrett}) and letting $w(x) = v(x) - \eta_{i}$
$$
J_{0}\left( (v)_{(\eta_{i}, \eta_{i+1})} \right) = \int_{\set{w(x) \in (0,\eta_{i+1}-\eta_{i})}}\ \left( \frac 1 2\, |\nabla w(x)|^{2} + R(w(x)+\eta_{i}) \right)\, dx \ \ge 0
$$
by Lemma \ref{esiste-kbar} since $R(s+\eta_{i})$ is positive for $s\to 0^{+}$. Since $\ell$ is finite, if $\| v - \tilde u \|$ is small enough then this is true at the same time for all $i= 1, \dots, \ell$.

Moreover, since $\tilde \omega < \Omega$, it follows that $K(\tilde u) > \frac{\sigma}{\Omega}$, hence the term
$$
\frac 1 2\, \left( \Omega^{2}\, K(u) + \frac{\sigma^{2}}{K(u)} \right)
$$
is strictly increasing in $K(\tilde u)$.

Putting all together, it follows that the energy $E_{\sigma}$ with nonlinear term $R$ satisfies for functions $v$ close to $\tilde u$ in the $H^{1}$ norm
$$
E_{\sigma}(v) = J_{0}\left( v_{(0,\eta_{1})} \right) + \sum_{i}\, J_{0}\left( (v)_{(\eta_{i}, \eta_{i+1})} \right) + \frac 1 2\, \left( \Omega^{2}\, K(v) + \frac{\sigma^{2}}{K(v)} \right) \ \ge
$$
$$
\ge \ J_{0}\left( v_{(0,\eta_{1})} \right) + \frac 1 2\, \left( \Omega^{2}\, K\left(v_{(0,\eta_{1})}\right) + \frac{\sigma^{2}}{K\left(v_{(0,\eta_{1})}\right)} \right) \ = \ \tilde E_{\sigma}(v) > \tilde E_{\sigma}(\tilde u) = E_{\sigma}(\tilde u)
$$
where we have used $\| \tilde u(x) \|_{L^{\infty}} < \eta_{1}$ and the fact that $\tilde E_{\sigma}$ and $E_{\sigma}$ coincide for functions with $L^{\infty}$ norm less than $\eta_{1}$ by (\ref{tilde-r}). \qed

\begin{lemma} \label{secondo-minimo}
There exists a threshold $\sigma_{g,2}$ such that, for all $\sigma > \sigma_{g,2}$, there exists a point of local minimum for the energy $E_{\sigma}(u)$ with $\xi_{2} < \| u(x) \|_{L^{\infty}} < \eta_{2}$.
\end{lemma}

\noindent \emph{Proof.}  Let us consider the modified nonlinear term
$$
\tilde R(s) = \left\{
\begin{array}{ll}
R(s) & s\le \eta_{2} \\[0.2cm]
f(s) & s\ge \eta_{2}
\end{array} \right.
$$
for any function $f(s)$ for which $f'(s) \ge 0$ and such that $\tilde R$ is of class $C^{2}$ and satisfies (\ref{H0})-(\ref{H3}). Then $\tilde R$ satisfies the assumptions of Lemma \ref{principio-massimo} with $\bar s = \eta_{2}$. As in the proof of Lemma \ref{primo-minimo}, we find a point of local minimum $u$ for the energy $\tilde E_{\sigma}$, and by Lemma \ref{principio-massimo} all points of local minimum satisfy $\| u(x) \|_{L^{\infty}} < \eta_{2}$. Then it follows that $u$ is a point of local minimum also for $E_{\sigma}$ and it satisfies $\| u(x) \|_{L^{\infty}} < \eta_{2}$.

It remains to show that there exists a point of local minimum for $\tilde E_{\sigma}$ with $
\| u(x) \|_{L^{\infty}} > \xi_{2}$. In the following we always assume to work with $\tilde R$ in the definition of $J_{0}$. Let us consider the set
$$
J_{(\eta_1,\eta_2)} := \jo \cap \set{J_0 \left( u_{(\eta_1,\eta_2)} \right) <0}
$$
which is open by Lemma \ref{cont-jo-ristr}. We first need to show that for $\sigma$ large enough
\begin{equation} \label{inf-interno-eta-2}
\inf\limits_{J_{(\eta_1,\eta_2)}}\ \tilde E_{\sigma}(u) < \inf\limits_{\partial J_{(\eta_1,\eta_2)}}\ \tilde E_{\sigma}(u)
\end{equation}
By Proposition \ref{prop:thresh}, we know that for all $\sigma$ large enough the (\ref{NKG}) equation admits hylomorphic solitons with charge $\sigma$. Hence, the infimum of $E_{\sigma}$ is not reached on functions $v$ on the boundary of $J_{(\eta_1,\eta_2)}$ for which $J_{0}(v)=0$. Indeed for such functions $v$, we have $\tilde E_{\sigma}(v) \ge \sigma\, \Omega$.

Moreover, by Lemma \ref{esiste-kbar}, there exists $\bar k_{(\eta_1,\eta_2)}>0$ such that functions $u \in \set{J_0 \left( u_{(\eta_1,\eta_2)} \right) <0}$ satisfy
$$
K \left( u_{(\eta_1,\eta_2)} \right) \ge \bar k_{(\eta_1,\eta_2)}
$$
hence by
$$
m\left( \set{u(x) \in (\eta_1,\eta_2)} \right)\, \eta_{2}^{2} \ge K \left( u_{(\eta_1,\eta_2)} \right)
$$
and Lemma \ref{conv-fort-caratt}, it follows that $v \in \partial J_{(\eta_1,\eta_2)}$ implies $\| v \|_{L^{\infty}} > \eta_{1}$. Hence, for $\sigma$ large enough, the infimum of $\tilde E_{\sigma}$ is reached neither on functions $v$ on the boundary of $J_{(\eta_1,\eta_2)}$ for which $J_{0}(v_{(\eta_1,\eta_2)})=0$. Indeed for such functions, letting $w(x) = v(x) - \eta_{1}$, it follows
$$
0= \int_{\set{w(x) \in (0,\eta_{2}-\eta_{1})}}\ \left( \frac 1 2\, |\nabla w(x)|^{2} + \tilde R(w(x)+\eta_{1}) \right)\, dx \le \int_{\R^{N}} \ \left( \frac 1 2\, |\nabla w|^{2} + \tilde R(w+\eta_{1}) \right)
$$
since $\tilde R(s+ \eta_{1})$ is non-negative for $s\ge \eta_{2}-\eta_{1}$. Hence again by 
Proposition \ref{prop:thresh}, we know that for all $\sigma$ large enough there exist hylomorphic solitons for the (\ref{NKG}) with nonlinear term $\tilde R(s + \eta_{1})$ which are in the set $\set{J_0 \left( u_{(0,\eta_2-\eta_{1})} \right) <0}$. Perturbing these solutions, by a cut-off on the tails, we find functions $u \in \set{J_0 \left( u_{(\eta_1,\eta_2)} \right) <0}$ with energy $\tilde E_{\sigma}$ strictly smaller than $\tilde E_{\sigma}(v)$. This proves (\ref{inf-interno-eta-2}).

Now, using (\ref{inf-interno-eta-2}) and Lemma \ref{p-s}, we obtain a point of local minimum for $\tilde E_{\sigma}$ in the open set $J_{(\eta_1,\eta_2)}$. That this point satisfies $\| u(x) \|_{L^{\infty}} > \xi_{2}$ is immediate from the definition of $J_{(\eta_1,\eta_2)}$ and $\xi_{2}$ (see (\ref{estr-int})). \qed

\vskip 0.5cm

The proof of Theorem \ref{main-3} for $\ell$ finite easily follows by repeating Lemma \ref{secondo-minimo} for all intervals $C_{i}$. Some modifications are needed for the last interval $C_{\ell}$ in the case $\eta_{\ell}= \infty$. There is no need for using a modified nonlinear term, since we only need to show the existence of a point of local minimum with $\| u(x) \|_{L^{\infty}} > \xi_{\ell}$. This is achieved as above by slightly modifying the argument. Finally, the interval of charges is given by
$$
\Sigma = \left( \max_{i=1,\dots,\ell}\, \sigma_{g,i},\, +\infty \right)
$$
If $\ell$ is infinite, we choose for any integer $M>0$ a modified nonlinear term with only $M$ intervals $C_{i}$ where it is negative, and such that it is positive for $s$ big enough. Then the result follows by repeating the same argument as above, showing in particular that the points of local minimum found for the modified energy are points of local minimum also for the original energy (follow the proof of Lemma \ref{primo-minimo}).



\begin{thebibliography}{99}

\bibitem{BBBM} J.Bellazzini, V.Benci, C.Bonanno, A.M.Micheletti, \emph{
Solitons for the nonlinear Klein-Gordon equation}, Adv. Nonlinear Stud., in press

\bibitem{hylo} J.Bellazzini, V.Benci, C.Bonanno, E.Sinibaldi, \emph{
Hylomorphic solitons in the nonlinear Klein-Gordon equations},
arXiv:0810.5079v1 [math.AP]

\bibitem{benci} V.Benci, G.Cerami, {\it The effect of the domain
topology on the number of positive solutions of nonlinear elliptic
problems}, Arch. Ration. Mech. Anal. \textbf{114} (1991), 79--93

\bibitem{manifold} V.Benci, C.Bonanno, A.M.Micheletti, \emph{On the multiplicity of solutions of a nonlinear elliptic problem on Riemannian manifolds}, J. Funct. Anal. \textbf{252} (2007), 464--489

\bibitem{Beres-Lions} H.Berestycki, P.L.Lions, \textit{Nonlinear scalar
field equations. I. Existence of a ground state}, Arch. Ration.
Mech. Anal. \textbf{82} (1982), 313--345

\bibitem{prep} C.Bonanno, \emph{Multiplicity of positive solutions for nonlinear field equations in $\R^{N}$}, arXiv:0901.4459v1 [math.AP]

\bibitem{coleman78} S.Coleman, V.Glaser, A.Martin, \emph{Action minima
among solutions to a class of euclidean scalar field equation},
Comm. Math. Phys. \textbf{58} (1978), 211--221

\bibitem{dancer} E.N.Dancer, {\it The effect of domain shape on
the number of positive solutions of certain nonlinear equations},
J. Differential Equations {\bf 74} (1988), 120--156

\bibitem{gss87} M.Grillakis, J.Shatah, W.Strauss, \emph{Stability theory of
solitary waves in the presence of symmetry, I}, J. Funct. Anal.
\textbf{74} (1987), 160--197

\bibitem{hirano} N.Hirano, \emph{Multiple existence of solutions for a nonlinear elliptic problem on a Riemannian manifold}, Nonlinear Anal. \textbf{70} (2008), 671--692

\bibitem{palais} R.S.Palais, \emph{The principle of symmetric criticality}, Comm. Math. Phys. \textbf{79} (1979), 19--30

\bibitem{rosen68}  G. Rosen, \emph{Particle-like solutions to nonlinear
complex scalar field theories with positive-definite energy densities}, J.
Math. Phys. \textbf{9} (1968), 996--998

\bibitem{shatah} J.Shatah, \emph{Stable standing waves of non-linear
Klein-Gordon equations}, Comm. Math. Phys. \textbf{91} (1983),
313--327

\bibitem{shatah2} J.Statah, \emph{Unstable ground state of nonlinear Klein-Gordon equation}, Trans. Amer. Math. Soc. \textbf{290} (1985), 701--710

\bibitem{ss85} J.Shatah, W.Strauss, \emph{Instability of nonlinear bound
states}, Comm. Math. Phys. \textbf{100} (1985), 173--190

\bibitem{strauss} W.A.Strauss, \emph{Nonlinear invariant wave equations},
Lecture notes in physics, vol. 23, Springer, 1978

\bibitem{VP} L.R. Volevic, B.P. Paneyakh, \emph{Certain spaces of generalized functions and embedding theorems}, Russian Math. Surv. \textbf{20} (1965), 1--73

\end{thebibliography}
\end{document}